\documentclass[a4paper,11pt,reqno]{amsart}
\usepackage{mathrsfs}
 \usepackage[all]{xy}
\usepackage{amsmath,amssymb,amscd,bbm,amsthm,mathrsfs}
\usepackage{color,xcolor}
\usepackage{graphicx}
\usepackage{manfnt}
\newtheorem{thm}{Theorem}[section]
\newtheorem{lem}{Lemma}[section]

\newtheorem{con}{Conjecture}[section]

\newtheorem{rem}{Remark}[section]

\begin{document}
\numberwithin{equation}{section}

 \title[Constant $k$th-mixed curvature]{Constant $k$th-mixed curvature}
\author{Weiguo Chen,\,\,\,\,Kai Tang}

\address{Weiguo Chen. School of Mathematical Sciences, Zhejiang Normal University, Jinhua, Zhejiang, 321004, China} \email{{chenwg@zjnu.edu.cn}}
\keywords{Hermitian manifold; Holomorphic sectional curvature; Chern Ricci curvature}
\address{Kai Tang. School of Mathematical Sciences, Zhejiang Normal University, Jinhua, Zhejiang, 321004, China} \email[Corresponding author]{{kaitang001@zjnu.edu.cn}}
\keywords{Hermitian manifold; Holomorphic sectional curvature; Chern Ricci curvature}
\thanks{\text{Foundation item:} Supported by National Natural Science
Foundation of China (Grant No.12001490).}
\begin{abstract}
In this paper, we consider general {\em kth-mixed curvature} $\mathcal{C}^{(k)}_{\alpha,\beta}$ ($\beta\neq0$) for Hermitian manifolds, which is  a convex combination of the $k$th Chern Ricci curvature and holomorphic sectional curvature. We prove that any compact Hermitian surface with constant $k$th-mixed curvature is self-dual. Furthermore, we show that if a compact Hermitian surface has constant
2th-mixed curvature $c$, then the Hermitian metric must be K\"{a}hler. For the higher-dimensional case, when the parameters $\alpha$ and $\beta$ satisfy certain conditions, we can also obtain partial results.
\end{abstract}

 \maketitle

\section{Introduction}
A complete K\"{a}hler manifold with constant holomorphic sectional curvature
 is called a complex space form. As we all know, they are quotients of complex projective space
 $\mathbb{CP}^{n}$, the complex Euclidean space $\mathbb{C}^{n}$, and the complex hyperbolic space $\mathbb{CH}^{n}$.

For a Hermitian manifold $(M^{n},g)$, the first Ricci curvature $Ric$  of the Chern connection is defined by
\begin{align}
Ric=Ric^{(1)}=\sqrt{-1}(g^{k\overline{l}}R_{i\overline{j}k\overline{l}})dz^{i}\wedge d\overline{z}^{j}=-\sqrt{-1}\partial\overline{\partial}\log \det g\,. \nonumber
\end{align}
where it is a $(1,1)$-form representing the first Chern class $c_{1}(M)$. The holomorphic sectional curvature $H$ is defined by
\begin{align}
H(X)=R(X,\overline{X},X,\overline{X})/|X|^{4} \nonumber \,\,,
\end{align}
for a (1, 0)-tangent vector $X\in T^{1,0}M$ (\cite{Zheng}). When the metric $g$ is non-K\"{a}hler, $H$ can not determine the entire curvature $R$ in general. A long term folklore conjecture in non-K\"{a}hler geometry is the following:
\begin{con}
	Let $(M^{n},g)$ be a compact Hermitian manifold with $n\geq2$. Assume that $H=c$ where $c$ is a constant. If $c\neq0$, then $g$ is K\"{a}hler and if $c=0$, then $R=0$.
\end{con}
Note that compact Chern flat manifolds were classified by Boothby \cite{Boothby} in 1958 as the set of
all compact quotients of complex Lie groups. In complex dimension $n=2$, the  conjecture holds. In the $c\leq0$
case, it was solved by Balas-Gauduchon \cite{BG} (see also \cite{Balas}) in 1985 for Chern connection and by Sato-Sekigawa \cite{SS} in 1990 for Riemannian connection. The general case for $n=2$ (for both Chern connection and
Levi-Civita connection) was solved by Apostolov-Davidov-Mu\v{s}karov \cite{ADM} in 1996, as a corollary of their beautiful
classification theorem for compact self-dual Hermitian surfaces.

For $n\geq3$, the conjecture is still largely open. In \cite{Tang1}, the second named
 author confirmed the conjecture under the additional assumption that
$g$ is Chern K\"{a}hler-like (namely $R$ obeys all K\"{a}hler symmetries). In \cite{CCN}, Chen-Chen-Nie proved the conjecture under the additional assumption that $g$ is locally conformally K\"{a}hler and $c\leq0$.  They also pointed out
 the necessity of the compactness assumption in the conjecture by explicit examples.  In \cite{ZhouZheng}, Zhou-Zheng  proved that any compact Hermitian threefold with vanishing
 {\em real bisectional curvature} must be Chern flat. Real bisectional curvature is a curvature notion
 introduced by  Yang-Zheng in \cite{YZ}. It is equivalent to holomorphic
 sectional curvature $H$ in strength when the metric is K\"{a}hler, but is slightly stronger than $H$
 when the metric is non-K\"{a}hler.  In \cite{RZ}, Rao-Zheng showed that the
 conjecture holds if the metric is Bismut K\"{a}hler-like, meaning that the curvature of the Bismut
 connection obeys all K\"{a}hler symmetries. Subsequently, Zhao-Zheng \cite{ZhaoZheng} proved a very beautiful theorem,  which states that all compact Strominger K\"{a}hler-like manifolds are pluriclosed. In \cite{BT3}, Broder and the second named
 author proved that pluriclosed metrics with vanishing holomorphic curvature on compact K\"{a}hler manifolds are K\"{a}hler. They also showed that Hermitian metrics with vanishing real bisectional curvature on compact complex manifolds in the Fujiki class $\mathcal{C}$ are K\"{a}hler.

Given the perplexing relationship between the holomorphic sectional curvature and
the first Ricci curvature,  one can attempt to interpolate between these curvatures
by considering the following curvature constraint: In \cite{CLT} (see also \cite{BT}), Chu-Lee-Tam
introduced, for $\alpha,\beta\in\mathbb{R}$ the {\em mixed curvature}
 \begin{align}
\mathcal{C}_{\alpha,\beta}(X)=\frac{\alpha}{|X|^{2}_{g}}Ric(X,\bar{X})+\beta H(X)\nonumber \,\,.
 \end{align}
We always assume that $\alpha,\beta$ are both not $0$, otherwise, it makes no sense. $\mathcal{C}_{1,0}$ is the Ricci curvature, $\mathcal{C}_{0,1}$ is the holomorphic sectional curvature. When metric $g$ is K\"{a}hler, mixed curvature encompasses many interesting curvature conditions. $\mathcal{C}_{1,1}$ is the notion $Ric^{+}(X,\overline{X})$ introduced by Ni \cite{N2}. $\mathcal{C}_{1,-1}$ is the orthogonal Ricci curvature $Ric^{\perp}(X,\overline{X})$  introduced by Ni-Zheng \cite{NZ1}. $\mathcal{C}_{k-1,n-k}$ is closely related to the $k$-Ricci curvature $Ric_{k}$ introduced by Ni \cite{N1}.

Recently, there are many important results on compact K\"{a}hler manifolds with $\mathcal{C}_{\alpha,\beta}\geq 0$ or $\mathcal{C}_{\alpha,\beta}\leq 0$. It was proved by Yang \cite{Yang2018} that a compact K\"{a}hler manifold with $\mathcal{C}_{0,1}>0$ must be projective and rationally connected, confirming a conjecture of Yau \cite{Yau}. In \cite{N2}, Ni showed that it is also true if $Ric_{k}>0$ for some $1\leq k\leq n$. In \cite{Mat1,Mat2}, Matsumura established the structure theorems for a projective manifold with $\mathcal{C}_{0,1}\geq0$. In \cite{CLT, BT},
it was shown that any compact K\"{a}hler manifold with $\mathcal{C}_{\alpha,\beta}>0$, which constants satisfy $\alpha>0 $ and $3\alpha+2\beta\geq0$, must be projective and simply connected.  Zhang-Zhang \cite{ZZ} generalized Yang's result to quasi-positive
 case, confirming a conjecture of Yang affirmatively.  The second named
 author \cite{Tang2} showed the projectivity of compact K\"{a}hler manifolds with quasi-positive $\mathcal{C}_{\alpha,\beta}$. Chu-Lee-Zhu \cite{CLZ} established a structure theorem for compact K\"{a}hler manifolds with $\mathcal{C}_{\alpha,\beta}\geq0$. Very recently, Zhang-Zhang\cite{ZZ1} proved that a compact K\"{a}hler manifold with quasi-positive $\mathcal{C}_{1,-1}$  is rationally connected, confirming a conjecture by Ni\cite{N2}.
 For the non-positive case, Wu-Yau \cite{WuYau} confirmed a conjecture of Yau that a projective K\"{a}hler manifold with $\mathcal{C}_{0,1}<0$ must have ample canonical line bundle. Tosatti-Yang \cite{TosattiYang} was able to drop the projectivity assumption in Wu-Yau theorem. Chu-Lee-Tam \cite{CLT} proved that a compact K\"{a}hler manifold with $\mathcal{C}_{\alpha,\beta}<0$ and $\alpha\geq0, \beta\geq 0$, have ample canonical bundle. For more related works, we refer readers to \cite{BT2,DT,HW,LS,N3,NZ2,Tang12,WuYau1,Yang20201,Zhang}.

Previous results have primarily focused on the study of K\"{a}hler cases, and much less is known about the non-K\"{a}hler case. In \cite{Tang3}, the second named author proposed the following conjecture:
\begin{con}\label{con2}
	Let $(M^{n},g)$ be a compact Hermitian manifold with $n\geq2$. Assume that mixed curvature $\mathcal{C}_{\alpha,\beta}=c$ where $c$ is a constant. If $c\neq0$, then $g$ is K\"{a}hler.
\end{con}
If Chern Ricci curvature is non-zero constant $\lambda$, namely $Ric=\lambda g$, then $g$ is  automatically K\"{a}hler. Note that if $\mathcal{C}_{\alpha,\beta}=0$, we may not be able to conclude that $R=0$. For example, isosceles Hopf manifold satisfies $\mathcal{C}_{\alpha,\beta}=0$ with $n\alpha+\beta=0$, but it does not imply $R=0$.

In complex dimension $n = 2$, the second named author\cite{Tang3} proved the Conjecture \ref{con2}. For $n\geq3$,
he also obtained results for locally conformally K\"{a}hler manifolds and for Chern K\"{a}hler-like manifolds. Chen-Zheng \cite{CZ1} verified the Conjecture \ref{con2} for several special types of Hermitian manifolds, including
 complex nilmanifolds, solvmanifolds with complex commutators, almost abelian Lie groups,
 and Lie algebras containing a $J$-invariant abelian ideal of codimension 2.

In the non-K\"{a}hler case, due to the lack of good symmetries in the curvature tensor, there are four distinct traces of the Chern curvature tensor,
thus, four distinct Chern Ricci curvatures. We denote by $Ric^{(k)}$ the $k$th Ricci curvature,  where the specific expression is given in Section 2. Then we define {\em kth-mixed curvature}:
 \begin{align}
\mathcal{C}^{(k)}_{\alpha,\beta}(X)=\frac{\alpha}{|X|^{2}_{g}}Ric^{(k)}(X,\bar{X})+\beta H(X)\nonumber \,\,,
 \end{align}
where $k=1,2,3,4$. In particular, when $k=1$, it is the mixed curvature. In this paper, we assume that $\beta\neq0$. Mimicking the above Conjecture \ref{con2}, we propose the following:
\begin{con}
	Let $(M^{n},g)$ be a compact Hermitian manifold with $n\geq2$. Assume that mixed curvature $\mathcal{C}^{(k)}_{\alpha,\beta}=c$ where $c$ is a constant. If $c\neq0$, then $g$ is K\"{a}hler.
\end{con}

In this paper, we first prove
\begin{thm}\label{th1.1}
Any compact Hermitian surface with constant $k$th-mixed curvature is self-dual.
\end{thm}

Apostolov-Davidov-Mu\v{s}karov \cite{ADM} proved that any compact self-dual surface is conformally to a special surface, which is either (1) a complex space form, or (2) the non-flat, conformally flat K\"{a}hler metric, or (3) an isosceles Hopf surface.

In particular, we also obtain
\begin{thm}\label{th1.2} If $(M^{2},g)$ is a compact Hermitian surface with constant 2th-mixed curvature $\mathcal{C}^{(2)}_{\alpha,\beta}=c$, then $g$ must be K\"{a}hler.
\end{thm}

\begin{rem} When $k=1$, the second named author \cite{Tang3} proved that if a compact Hermitian surface $(M,g)$ has constant 1th-mixed curvature $\mathcal{C}^{(1)}_{\alpha,\beta}=c$, then $g$ must be K\"{a}hler unless $c=0$ and $2\alpha+\beta=0$. Isosceles Hopf surface satisfies $\mathcal{C}^{(1)}_{\alpha,\beta}=0$ with $2\alpha+\beta=0$, where it is non-K\"{a}hler manifold. However, the above Theorem \ref{th1.2} holds for arbitrary constants $\alpha$ and $\beta$. Note that if $\alpha\neq0$ and $\beta=0$, we have $Ric^{(2)}=cg$, where the isosceles Hopf surface satisfies this condition (see section 2).
\end{rem}

For higher-dimensional cases ($n\geq3$), we can discuss some special cases.
\begin{thm}\label{th1.3}
Let $(M^{n},g)$ be a compact Hermitian manifold with $\mathcal{C}^{(k)}_{\alpha,\beta}=0$. If $\alpha(n+1)+2\beta=0$, then $g$ is balanced metric.
\end{thm}
Therefore, we have that any compact Hermitian surface with $\mathcal{C}^{(k)}_{2,-3}=0$ must be K\"{a}hler.

This paper is organized as follows. In section 2, we provide some basic knowledge which will be used in our proofs. In section 3, we prove Theorem \ref{th1.1}, \ref{th1.2}, \ref{th1.3}.

\section{Preminaries}
Let $(M,g)$ be a compact Hermitian manifold of dimension $n$ and denote by $\omega$ the K\"{a}hler form associated with $g$. If $\omega$ is closed, that is, if $d\omega=0$, we call $g$ a K\"{a}hler metric. Denote by $\nabla$ Chern connection. If we choose local holomorphic chart $(z_{1},\cdot\cdot\cdot,z^{n})$, we write
\begin{align}
\omega=\sqrt{-1}\sum_{i=1,j=1}^{n}g_{i\overline{j}}dz^{i}\wedge d\overline{z}^{j}. \nonumber
\end{align}
Then the curvature tensor $R=\{R_{i\overline{j}k\overline{l}}\}$ of the Chern connection is given by
\begin{align}
R_{i\overline{j}k\overline{l}}=-\frac{\partial^{2}g_{k\overline{l}}}{\partial z^{i}\partial \overline{z}^{j}}+g^{p\overline{q}}\frac{\partial g_{k\overline{q}}}{\partial z^{i}}\frac{\partial g_{p\overline{l}}}{\partial \overline{z}^{j}}   . \nonumber
\end{align}
There are four distinct traces of the Chern curvature tensor and, thus, four distinct Chern Ricci curvatures.  The \textit{first} and \textit{second Chern Ricci curvatures} $\text{Ric}_{i \bar{j}}^{(1)}  : =R_{i\overline{j}}=  \text{g}^{k \bar{\ell}} \text{R}_{i \bar{j} k \bar{\ell}}$ and  $\text{Ric}_{k \bar{\ell}}^{(2)}  : =R_{kl}^{(2)}=  \text{g}^{i \bar{j}} \text{R}_{i \bar{j} k \bar{\ell}}$ trace to the same \textit{Chern scalar curvature} $u$.  The remaining \textit{third} and \textit{fourth Chern Ricci curvatures} $\text{Ric}_{i \bar{\ell}}^{(3)}  : =  \text{g}^{k \bar{j}} \text{R}_{i \bar{j} k \bar{\ell}}$ and $\text{Ric}_{k \bar{j}}^{(4)} : =  \text{g}^{i \bar{\ell}} \text{R}_{i \bar{j} k \bar{\ell}}$ trace to what we refer to as the \textit{altered Chern scalar curvature} $v$, i.e.,   \begin{eqnarray*}
 u\ : = \ \text{g}^{i \bar{j}} \text{Ric}_{i \bar{j}}^{(1)} \ = \ \text{g}^{k \bar{\ell}} \text{Ric}_{k \bar{\ell}}^{(2)}, \qquad v\ : = \ \text{g}^{i \bar{\ell}} \text{Ric}_{i \bar{\ell}}^{(3)} \ = \ \text{g}^{k \bar{j}} \text{Ric}_{k \bar{j}}^{(4)}.
\end{eqnarray*}
The first Chern Ricci curvature is a $(1,1)$-form representing the first Chern class $c_{1}(M)$. For convenience, we also denote the first, second, third, and fourth Ricci curvatures as $\rho^{(1)}$, $\rho^{(2)}$, $\rho^{(3)}$ and $\rho^{(4)}$ respectively.
 The components of the Chern torsion T are given by $\text{T}_{ij}^k = \text{g}^{k \bar{\ell}} (\partial_i \text{g}_{j \bar{\ell}} - \partial_j \text{g}_{i \bar{\ell}})$.  We write $\eta = \sum_i \eta_i dz^i = \sum_{i,k} \text{T}_{ik}^k dz^i$ for the \textit{torsion $(1,0)$-form}.  A Hermitian metric is called balanced if $d\omega^{n-1}=0$, namely, $\eta=0$.

It is well know that
 \begin{align}
\int_{M}(u-v)d\nu=\int_{M}|\eta|^{2}d\nu \nonumber
\end{align}
 Note that if $g$ is K\"{a}hler metric, then $R$ obeys all K\"{a}hler symmetries, four Chern Ricci curvatures are same and T$=0$.

By the conformal formula for Hermitian metric, the Hermitian $\widetilde{g}=e^{2F}g$ satisfies
 \begin{align}\label{2.01}
\widetilde{R}_{k\overline{l}i\overline{j}}=e^{2F}(R_{k\overline{l}i\overline{j}}-2g_{i\overline{j}}F_{k\overline{l}})
\end{align}
where $\widetilde{R}$ is the curvature tensor of $\widetilde{g}$ and $R$ is the curvature tensor of $g$.

For Hermitian surfaces $(M^{2},g)$, we list a few special formulas. By \cite[p113,(4.2)]{KO}, we get
 \begin{align}\label{2.2}
c_{1}^{2}(M)=\frac{1}{8\pi^{2}}[u^{2}-\langle Ric^{(1)},Ric^{(1)}\rangle]g^{2}
\end{align}
From \cite{BG} (also see \cite{G}), we have
\begin{align}\label{2.3}
&\rho^{(1)}+\rho^{(2)}-2 Re(\rho^{(3)})=(u-v)g \\
&\int_{M}\langle \rho^{(1)},\rho^{(2)}\rangle g^{2}=\int_{M}\langle \rho^{(1)}, Re(\rho^{(3)})\rangle g^{2}=
\int_{M}\langle\rho^{(1)},\rho^{(1)}\rangle g^{2}-\int_{M}u(u-v)g^{2}
\end{align}

Next let us recall the formula of Weyl curvature tensor for a Hermitian surface $(M^{2},g)$. One has $W=W^{+}+W^{-}$, and $g$ is said to be {\em self-dual} if $W^{-}=0$. Let $\{e_{1},e_{2}\}$ be a local unitary frame of $M$. Then
\begin{align}
\{e_{1}\wedge \overline{e}_{2}, e_{2}\wedge \overline{e}_{1}, \frac{1}{\sqrt{2}}(e_{1}\wedge\overline{e}_{1}-e_{2}\wedge\overline{e}_{2})\} \nonumber
\end{align}
form a basis of the complexification of $\Lambda^{2}_{-}$, the space on which the Hodge star operator is minus identity. We  call this basis the standard basis associated with the unitary frame $e$. In \cite{ADM}, it was proved that
\vspace{0.2cm}
\begin{lem}[\cite{ADM}]\label{lemma2.1}
Let $(M^{2},g)$ be a Hermitian surface, and $e$ a local unitary frame. Then under the standard basis associated with $e$, the components of $W^{-}$, the anti-self dual part of the Weyl tensor, are given by
\begin{align}
&W_{1}^{-}=R_{1\overline{2}1\overline{2}    }                                      \nonumber \\
&W_{2}^{-}= \frac{1}{\sqrt{2}}(R_{1\overline{2}2\overline{2}}+R_{2\overline{2}1\overline{2}}-
R_{1\overline{2}1\overline{1}}-R_{1\overline{1}1\overline{2}})                                          \nonumber \\
&W_{3}^{-}=\frac{1}{6}(R_{1\overline{1}1\overline{1}}+R_{2\overline{2}2\overline{2}}-
R_{1\overline{1}2\overline{2}}-R_{2\overline{2}1\overline{1}}-R_{1\overline{2}2\overline{1}}-R_{2\overline{1}1\overline{2}})                                          \nonumber
\end{align}
\end{lem}

Recall that Hopf manifold are compact complex manifolds whose universal cover is $\mathbb{C}^{n}\setminus\{0\}$ (see
\cite{CZ}). An isosceles Hopf manifold means $M^{n}_{\phi}=\mathbb{C}^{n}\setminus\{0\}/\langle \phi\rangle$ where
\begin{align}
\phi: (z_{1},z_{2},\cdot\cdot\cdot,z_{n})\rightarrow (a_{1}z_{1},a_{2}z_{2},\cdot\cdot\cdot,a_{n}z_{n})
\end{align}
with $0<|a_{1}|=|a_{2}|=\cdot\cdot\cdot=|a_{n}|<1$. The standard Hopf metric $g$ with K\"{a}hler form $\omega_{g}=\sqrt{-1}\frac{\partial\overline{\partial}|z|^{2}}{|z|^{2}}$ descends down to $M^{n}_{\phi}$. Write $h=|z|$, $e_{i}=h\frac{\partial}{\partial z_{i}}$, $\varphi_{i}=\frac{1}{h}dz_{i}$. Then $e$ becomes a unitary frame with $\varphi$ the
dual coframe. We have
\begin{align}
R_{i\overline{j}k\overline{l}}=-\frac{\overline{z}_{i}z_{j}}{|z|^{2}}\delta_{kl}+\delta_{ij}\delta_{kl}
\end{align}
where $h_{i}=\frac{\partial h}{\partial z^{i}}$. From this, the first Ricci curvature is
\begin{align}\label{2.7}
R_{i\overline{j}}=-n\frac{\overline{z}_{i}z_{j}}{|z|^{2}}+n\delta_{ij}
\end{align}
and the second Ricci curvature is
\begin{align}\label{2.8}
R_{i\overline{j}}^{(2)}=(n-1)\delta_{i\overline{j}}
\end{align}
So we get scalar curvature $u=n^{2}-n$ and altered  scalar curvature $v=n-1$.
It implies
\begin{align}
(R_{i\bar{j}}g_{k\bar{l}}+R_{k\bar{j}}g_{i\bar{l}}+
R_{i\bar{l}}g_{k\bar{j}}+R_{k\bar{l}}g_{i\bar{j}})
-n(R_{i\bar{j}k\bar{l}}+R_{k\bar{j}i\bar{l}}+R_{i\bar{l}k\bar{j}}+R_{k\bar{l}i\bar{j}})
=0
\end{align}
So when $n\alpha+\beta=0$, the Hopf manifold satisies that $\mathcal{C}^{(1)}_{\alpha,\beta}=0$. By (\ref{2.8}), we also have that $Ric^{(2)}=(n-1)g$, where $g$ is the second Chern Einstein metrics. Note that the Hopf manifold is not Chern flat.

\section{Proofs of Theorem \ref{th1.1}, \ref{th1.2} and \ref{th1.3}}
\subsection{Hermitian surface} We will first prove Theorem \ref{th1.1}.

\noindent{\bf{{\em Proof of Theorem \ref{th1.1}.}}}
Let $(M^{2},g)$ be a  Hermitian surface. Assume that
 \begin{align}
\mathcal{C}^{(k)}_{\alpha,\beta}(X)=\frac{\alpha}{|X|^{2}_{g}}Ric^{(k)}(X,\bar{X})+\beta H(X)=c \nonumber \,\,,
 \end{align}
 where $c$ is constant and $k=1,2,3,4$. When $k=1$, this result has already been proved by the second named author \cite{Tang3}. We mainly consider the cases where $k=2, 3, 4$. For any fixed $p\in M$, let $X=e_{1}$ and $Y=e_{2}$ be orthonormal $(1,0)$-vector.

(1) The first case $k=2$. We have
\begin{align*}
&\frac{1}{vol(\mathbb{S}^{3})}\int_{ Z \in T^{1,0}_{p}M, | Z | =1 }\alpha Ric^{(k)}(Z,\overline{Z})+\beta H(Z) d\theta(Z) \nonumber \\
&=
\frac{\alpha}{2}u+\frac{\beta}{6}(u+v)=\frac{(3\alpha+\beta)u+\beta v}{6}=c \nonumber
\end{align*}
where the first equality is based on the average trick.
Hence
\begin{align}\label{3.1}
(3\alpha+2\beta)(R_{1\overline{1}1\overline{1}}+R_{2\overline{2}2\overline{2}})+(3\alpha+\beta)(R_{1\overline{1}2\overline{2}}+
R_{2\overline{2}1\overline{1}})+\beta(R_{1\overline{2}2\overline{1}}+R_{2\overline{1}1\overline{2}})=6c
\end{align}
By $\mathcal{C}^{(2)}_{\alpha,\beta}(X)=\mathcal{C}^{(2)}_{\alpha,\beta}(Y)=c$, we have
\begin{align}\label{3.2}
&\alpha(R_{1\overline{1}1\overline{1}}+R_{2\overline{2}1\overline{1}})+\beta R_{1\overline{1}1\overline{1}}=c \\
&\alpha(R_{1\overline{1}2\overline{2}}+R_{2\overline{2}2\overline{2}})+\beta R_{2\overline{2}2\overline{2}}=c
\end{align}
From (\ref{3.1}), (\ref{3.2}), (3.3) and $\beta\neq0$, we get
\begin{align}\label{3.4}
R_{1\overline{1}1\overline{1}}+R_{2\overline{2}2\overline{2}}-(R_{1\overline{1}2\overline{2}}+R_{2\overline{2}1\overline{1}}+
R_{1\overline{2}2\overline{1}}+R_{2\overline{1}1\overline{2}})=0
\end{align}
So $W_{3}^{-}=0$.

(2) The second case $k=3, 4$. The average trick gives that
\begin{align*}
&\frac{1}{vol(\mathbb{S}^{3})}\int_{ Z \in T^{1,0}_{p}M, | Z | =1 }\alpha Ric^{(k)}(Z,\overline{Z})+\beta H(Z) d\theta(Z) \nonumber \\
&=
\frac{\alpha}{2}v+\frac{\beta}{6}(u+v)=\frac{(3\alpha+\beta)v+\beta u}{6}=c \nonumber
\end{align*}
So we have
\begin{align}\label{3.5}
(3\alpha+2\beta)(R_{1\overline{1}1\overline{1}}+R_{2\overline{2}2\overline{2}})+(3\alpha+\beta)(R_{1\overline{2}2\overline{1}}+
R_{2\overline{1}1\overline{2}})+\beta(R_{1\overline{1}2\overline{2}}+R_{2\overline{2}1\overline{1}})=6c
\end{align}
From $\mathcal{C}^{(3)}_{\alpha,\beta}(X)=\mathcal{C}^{(3)}_{\alpha,\beta}(Y)=c$, we obtain
\begin{align}\label{3.6}
&\alpha(R_{1\overline{1}1\overline{1}}+R_{1\overline{2}2\overline{1}})+\beta R_{1\overline{1}1\overline{1}}=c \\
&\alpha(R_{2\overline{1}1\overline{2}}+R_{2\overline{2}2\overline{2}})+\beta R_{2\overline{2}2\overline{2}}=c
\end{align}
By (\ref{3.5}), (\ref{3.6}), (3.7) and $\beta\neq0$, we also have
\begin{align}\label{3.8}
R_{1\overline{1}1\overline{1}}+R_{2\overline{2}2\overline{2}}-(R_{1\overline{1}2\overline{2}}+R_{2\overline{2}1\overline{1}}+
R_{1\overline{2}2\overline{1}}+R_{2\overline{1}1\overline{2}})=0
\end{align}
Therefore, we can also obtain $W_{-}^{3}=0$. The case for $k=4$ is similar to that of $k=3$.

Next, we will further prove that $W_{1}^{-}=0$ and $W_{2}^{-}=0$.

Let $a,b\in\mathbb{C}$ and $|a|^{2}+|b|^{2}=1$. Then $aX+bY$ is unit vector. Thus, by carrying out the calculation, we obtain
\begin{align}\label{3.9}
c&=\alpha Ric^{(k)}(aX+bY,\overline{a}\overline{X}+\overline{b}\overline{Y})+\beta R(aX+bY,\overline{a}\overline{X}+\overline{b}\overline{Y},
aX+bY,\overline{a}\overline{X}+\overline{b}\overline{Y}) \nonumber\\
&=\alpha\{|a|^{2}Ric^{(k)}(X,\overline{X})+a\overline{b}Ric^{(k)}(X,\overline{Y})+b\overline{a}Ric^{(k)}(Y,\overline{X})+|b|^{2}Ric^{(k)}(Y,\overline{Y})\} \nonumber\\
&\,\,\,\,\,\,+\beta\{|a|^{4}R(X,\overline{X},X,\overline{X}+|b|^{4}R(Y,\overline{Y},Y,\overline{Y})) \nonumber \\
&\,\,\,\,\,\,+|a|^{2}|b|^{2}[R(X,\overline{X},Y,\overline{Y})+R(Y,\overline{Y},X,\overline{X})+R(X,\overline{Y},Y,\overline{X})+
R(Y,\overline{X},X,\overline{Y})]\nonumber \\
&\,\,\,\,\,\,+a^{2}\overline{b}^{2}R(X,\overline{Y},X,\overline{Y})+\overline{a}^{2}b^{2}R(Y,\overline{X},Y,\overline{X}) \nonumber\\
&\,\,\,\,\,\,+a^{2}\overline{a}\overline{b}[R(X,\overline{X},X,\overline{Y})+R(X,\overline{Y},X,\overline{X})]
+\overline{a}^{2}ab[R(Y,\overline{X},X,\overline{X})+R(X,\overline{X},Y,\overline{X})]   \nonumber \\
&\,\,\,\,\,\,+\overline{a}b^{2}\overline{b}[R(Y,\overline{Y},Y,\overline{X})+R(Y,\overline{X},Y,\overline{Y})]
+a\overline{b}^{2}b[R(Y,\overline{Y},X,\overline{Y})+R(X,\overline{Y},Y,\overline{Y})]  \}
\end{align}
In fact, by combining (\ref{3.4}), $\mathcal{C}^{(k)}_{\alpha,\beta}(X)=\mathcal{C}^{(k)}_{\alpha,\beta}(Y)=c$ and $|a|^{2}+|b|^{2}=1$, (\ref{3.9}) can be simplified to the following
\begin{align}
0=&\alpha[a\overline{b}Ric^{(k)}(X,\overline{Y})+b\overline{a}Ric^{(k)}(Y,\overline{X})]+\beta
[a^{2}\overline{b}^{2}R(X,\overline{Y},X,\overline{Y})+\overline{a}^{2}b^{2}R(Y,\overline{X},Y,\overline{X})] \nonumber \\
&+\beta[(a^{2}\overline{a}\overline{b})A+(\overline{a}^{2}ab)B+(\overline{a}b^{2}\overline{b})C+(a\overline{b}^{2}b)D]
\end{align}
where we assume that
\begin{align}
&A=R(X,\overline{X},X,\overline{Y})+R(X,\overline{Y},X,\overline{X});\,\,\,\,
B=R(Y,\overline{X},X,\overline{X})+R(X,\overline{X},Y,\overline{X}) \nonumber \\
&C=R(Y,\overline{Y},Y,\overline{X})+R(Y,\overline{X},Y,\overline{Y});\,\,\,\,\,\,\,\,
D=R(Y,\overline{Y},X,\overline{Y})+R(X,\overline{Y},Y,\overline{Y})
\end{align}
Now we choose special $a$ and $b$. Setting $a=\frac{1}{\sqrt{2}}, b=\frac{1}{\sqrt{2}}$, and $a=\frac{-1}{\sqrt{2}}, b=\frac{1}{\sqrt{2}}$. Then we get
\begin{align}
R(X,\overline{Y},X,\overline{Y})+R(Y,\overline{X},Y,\overline{X})=0
\end{align}
Setting $a=\frac{\sqrt{2}}{2}, b=\frac{1-i}{2}$, and $a=\frac{-\sqrt{2}}{2}, b=\frac{1-i}{2}$. We have
\begin{align}
R(X,\overline{Y},X,\overline{Y})-R(Y,\overline{X},Y,\overline{X})=0
\end{align}
Hence $R_{1\overline{2}1\overline{2}}=R(X,\overline{Y},X,Y)=0$. So $W_{1}^{-}=0$. Now  (\ref{3.9}) is further simplified to
\begin{align}
0=&\alpha[a\overline{b}Ric^{(k)}(X,\overline{Y})+b\overline{a}Ric^{(k)}(Y,\overline{X})]
+\beta[(a^{2}\overline{a}\overline{b})A+(\overline{a}^{2}ab)B+(\overline{a}b^{2}\overline{b})C+(a\overline{b}^{2}b)D]
\end{align}
Setting $a=\frac{\sqrt{2}}{2}, b=-\frac{\sqrt{2}}{2}i$, and $a=\frac{\sqrt{3}}{2}, b=\frac{1}{2}i$, yields
\begin{align}
A+C=B+D \nonumber
\end{align}
Setting $a=-\frac{\sqrt{2}}{2}, b=\frac{\sqrt{2}}{2}$, and $a=\frac{\sqrt{3}}{2}, b=\frac{1}{2}$, yields
\begin{align}
A+B=C+D \nonumber
\end{align}
From these, we get
\begin{align}
R_{1\overline{2}2\overline{2}}+R_{2\overline{2}1\overline{2}}-R_{1\overline{2}1\overline{1}}-R_{1\overline{1}1\overline{2}}=D-A=0 \nonumber
\end{align}
namely, $W_{2}^{-}=0$. We get that $W^{-}=0$. According to Lemma \ref{lemma2.1}, we conclude that Hermitian surface $M$ must be self-dual.
\qed

\vspace{0.2cm}
Before proving Theorem \ref{th1.2}, we first present several preliminary lemmas. Although the proof idea of Theorem \ref{th1.2} is similar to that of Theorem 1.2 in \cite{Tang3}, some details are different.
\begin{lem}
Let $(M^{n},g)$ be a compact Hermitian manifold. Assume that 2th-mixed curvature $\mathcal{C}^{(2)}_{\alpha,\beta}=c$. Then we have
\begin{align}\label{3.15}
[\alpha(n+2)+\beta]Ric^{(2)}+\beta Ric^{(1)}+2\beta Re (Ric^{(3)})=[2(n+1)c-\alpha u]g
\end{align}
\end{lem}

\noindent{{\em Proof.}} For $p\in M$, Let $X=X^{i}\frac{\partial}{\partial x^{i}}\in T_{p}^{(1,0)}M$. Since $\mathcal{C}^{(2)}_{\alpha,\beta}(X)=c$ at $p$, then
\begin{align}
\alpha Ric^{(2)}(X,\overline{X})|X|^{2}+\beta R(X,\overline{X},X,\overline{X})=c|X|^{4}
\end{align}
So we have
\begin{align}
&\alpha X^{i}\overline{X}^{j}X^{k}\overline{X}^{l}(R^{(2)}_{i\bar{j}}g_{k\bar{l}}+R^{(2)}_{k\bar{j}}g_{i\bar{l}}+
R^{(2)}_{i\bar{l}}g_{k\bar{j}}+R^{(2)}_{k\bar{l}}g_{i\bar{j}}) \nonumber \\
&+\beta X^{i}\overline{X}^{j}X^{k}\overline{X}^{l}(R_{i\bar{j}k\bar{l}}+R_{k\bar{j}i\bar{l}}+R_{i\bar{l}k\bar{j}}+R_{k\bar{l}i\bar{j}})
=2cX^{i}\overline{X}^{j}X^{k}\overline{X}^{l}(g_{i\bar{j}}g_{k\bar{l}}+g_{i\bar{l}}g_{k\bar{j}})
\end{align}
This is equivalent to saying that
\begin{align}\label{3.18}
\alpha (R^{(2)}_{i\bar{j}}g_{k\bar{l}}+R^{(2)}_{k\bar{j}}g_{i\bar{l}}+
R^{(2)}_{i\bar{l}}g_{k\bar{j}}+R^{(2)}_{k\bar{l}}g_{i\bar{j}})
+\beta(R_{i\bar{j}k\bar{l}}+R_{k\bar{j}i\bar{l}}+R_{i\bar{l}k\bar{j}}+R_{k\bar{l}i\bar{j}})
=2c(g_{i\bar{j}}g_{k\bar{l}}+g_{i\bar{l}}g_{k\bar{j}})
\end{align}
Taking the trace,  and combining with $Ric^{(3)}+Ric^{(4)}=2Re (Ric^{(3)})$, we obtain (\ref{3.15}).                \qed

We will next consider the case of conformal change.
\begin{lem}\label{lem3.2}
Let $(M^{n},g)$ be a compact Hermitian manifold. Let $\widetilde{g}=e^{2F}g$ is a conformal metric with
constant 2th-mixed curvature $\mathcal{C}^{(2)}_{\alpha,\beta}=c$. We denote  curvatures of $g$ as $R_{i\overline{j}k\overline{l}}$. Then we have
\begin{align}\label{3.37}
&\alpha (R^{(2)}_{i\bar{j}}g_{k\bar{l}}+R^{(2)}_{k\bar{j}}g_{i\bar{l}}+
R^{(2)}_{i\bar{l}}g_{k\bar{j}}+R^{(2)}_{k\bar{l}}g_{i\bar{j}})
+\beta(R_{i\bar{j}k\bar{l}}+R_{k\bar{j}i\bar{l}}+R_{i\bar{l}k\bar{j}}+R_{k\bar{l}i\bar{j}}) \nonumber \\
&-2\beta(g_{i\overline{j}}F_{k\overline{l}}+g_{k\overline{j}}F_{i\overline{l}}+g_{i\overline{l}}F_{k\overline{j}}
+g_{k\overline{l}}F_{i\overline{j}})
=(4\alpha \Delta F+2c e^{2F})(g_{i\bar{j}}g_{k\bar{l}}+g_{i\bar{l}}g_{k\bar{j}})
\end{align}
\end{lem}

\noindent{\em Proof.} From formula (\ref{2.01}) and (\ref{3.18}), we can easily draw the conclusion. \qed

\begin{lem}\label{lem3.3}
Let $(M^{2},g)$ be a compact K\"{a}hler surface of constant holomorphic sectional curvature $H=c_{0}$. Let $\widetilde{g}=e^{2F}g$ is a conformal metric with with constant 2th-mixed curvature $\mathcal{C}^{(2)}_{\alpha,\beta}=c$. Then $F$ is constant. In particular, $\widetilde{g}$ must be a K\"{a}hler metric.
\end{lem}

\noindent{{\em Proof.}} We denote the first, second  and third Ricci curvatures of $\widetilde{g}$ as $\widetilde{\rho}^{(1)},\widetilde{\rho}^{(2)},\widetilde{\rho}^{(3)}$  respectively. Let $u,v$ is scalar and altered  scalar curvature of $g$, and $\widetilde{u},\widetilde{v}$ is scalar and altered  scalar curvature of $\widetilde{g}$.

By (\ref{2.3}) and (\ref{3.15}), we have
\begin{align}
&\widetilde{\rho}^{(1)}+\widetilde{\rho}^{(2)}-2 Re (\widetilde{\rho}^{(3)})=(\widetilde{u}-\widetilde{v})\widetilde{g}  \\
&[4\alpha+\beta]\widetilde{\rho}^{(2)}+\beta\widetilde{\rho}^{(1)}+2\beta Re(\widetilde{\rho}^{(3)})=[6c-\alpha\widetilde{u}]\widetilde{g}
\end{align}
and  since $[3\alpha+\beta]\widetilde{u}+\beta\widetilde{v}=6c$, we obtain
\begin{align}\label{3.22}
2\alpha\widetilde{\rho}^{(2)}+2\beta Re (\widetilde{\rho}^{(3)})=(\alpha\widetilde{u}+\beta\widetilde{v})\widetilde{g}
\end{align}
By (2.4),
\begin{align}\label{3.23}
\int_{M}\langle\widetilde{\rho}^{(1)},\widetilde{\rho}^{(2)} \rangle\widetilde{g}^{2}=
\int_{M}\langle\widetilde{\rho}^{(1)},Re (\widetilde{\rho}^{(3)})\rangle \widetilde{g}^{2}=\int_{M}\langle\widetilde{\rho}^{(1)},
\widetilde{\rho}^{(1)} \rangle\widetilde{g}^{2}-\int_{M}\widetilde{u}(\widetilde{u}-\widetilde{v})\widetilde{g}^{2}
\end{align}
Combining (\ref{3.22}) and (\ref{3.23}), we get
\begin{align}
(3\alpha+2\beta)\int_{M}\widetilde{u}^{2}\widetilde{g}^{2}=2(\alpha+\beta)\int_{M}\langle\widetilde{\rho}^{(1)},
\widetilde{\rho}^{(1)} \rangle\widetilde{g}^{2}+(2\alpha+\beta)\int_{M}(\widetilde{u}\widetilde{v})\widetilde{g}^{2}
\end{align}
Setting $c_{1}^{2}=\int_{M}c_{1}^{2}(M)$, according to (\ref{2.2}), we have
\begin{align}\label{3.25}
c_{1}^{2}=\frac{1}{8\pi^{2}}[\int_{M}\widetilde{u}^{2}\widetilde{g}^{2}-\int_{M}\langle\widetilde{\rho}^{(1)},
\widetilde{\rho}^{(1)} \rangle\widetilde{g}^{2}]
\end{align}
Hence, we have
\begin{align}\label{3.26}
16\pi^{2}(\alpha+\beta)c_{1}^{2}=-\alpha\int_{M}\widetilde{u}^{2}\widetilde{g}^{2}+(2\alpha+\beta)\int_{M}(\widetilde{u}\widetilde{v})\widetilde{g}^{2}
\end{align}
Since $g$ is K\"{a}hler metric with $H=c_{0}$, then $u=v=3c_{0}$. So we obtain
\begin{align}\label{3.27}
c_{1}^{2}=\frac{1}{16\pi^{2}}\int_{M}u^{2}g^{2}
\end{align}
From (\ref{3.26}) and (\ref{3.27}), we have
\begin{align}
(\alpha+\beta)\int_{M}u^{2}g^{2}=-\alpha\int_{M}\widetilde{u}^{2}\widetilde{g}^{2}+(2\alpha+\beta)\int_{M}(\widetilde{u}\widetilde{v})\widetilde{g}^{2}
\end{align}
Since $\widetilde{g}=e^{2F}g$, then by (\ref{2.01}), we obtain
\begin{align}\label{3.29}
e^{2F}\widetilde{u}=u-4\Delta F;\,\,\,\,\,\,e^{2F}\widetilde{v}=v-2\Delta F
\end{align}
From these, we get
\begin{align}
(\alpha+\beta)\int_{M}u^{2}g^{2}=-\alpha\int_{M}(u-4\Delta F)^{2}g^{2}+(2\alpha+\beta) \int_{M}(u-4\Delta F)(u-2\Delta F)g^{2}
\end{align}
Then
\begin{align}
\beta\int_{M}(\Delta F)^{2}g^{2}=0
\end{align}
where we use $\int_{M}\Delta F g^{2}=0$ and $u=3c_{0}$. Since $\beta\neq0$, then $\Delta F=0$. So $F$ is constant. \qed

\vspace{0.3cm}

\begin{lem}\label{lem3.4}
Let $M^{2}$ be a Hopf surface with standard Hopf metric $g$. Let $\widetilde{g}=e^{2F}g$ is a conformal metric with constant 2th-mixed curvature $\mathcal{C}^{(2)}_{\alpha,\beta}=c$. Then  metric $\widetilde{g}$ does not exist.
\end{lem}
\noindent{{\em Proof.}} From (\ref{2.2}) and (\ref{3.25}), we know that Hopf surface $c_{1}^{2} = 0$. By (\ref{3.26}), (\ref{3.29}) and $u=2v=1$, we have $\widetilde{u}=2\widetilde{v}$ and
\begin{align}
0=-4\alpha\int_{M}(v-2\Delta F)^{2}g^{2}+2(2\alpha+\beta)\int_{M}(v-2\Delta F)^{2}g^{2}=2\beta\int_{M}(v-2\Delta F)^{2}g^{2}
\end{align}
since $\beta\neq0$, then $\Delta F=\frac{v}{2}=\frac{1}{2}$, namely, $F$ is constant.  This of course will contradict with $v\neq0$. \qed

\vspace{0.2cm}
Now we are ready to prove the Theorem \ref{th1.2}.

\vspace{0.2cm}
\noindent{\bf{{\em Proof of Theorem \ref{th1.2}.}}}
Let $(M^{2},g)$ is a compact Hermitian surface with constant 2th-mixed curvature $\mathcal{C}^{(2)}_{\alpha,\beta}=c$. By Theorem \ref{th1.1}, we know that
$M^{2}$ is is self-dual. Hence, by Theorem 1' in \cite{ADM}, $g$ is conformally to a metric $h$ on $M^{2}$, which is either (1) a complex space form, or (2) the non-flat, conformally flat K\"{a}hler metric, or (3) an isosceles Hopf surface. Write $g=e^{2F}h$.

For case (1), by Lemma \ref{lem3.3}, we get that $F$ is constant. Hence $g$ is K\"{a}hler.

For case (2), $M^{2}$ is a holomorphic bundle over a curve $C$ of genus at least 2, and the K\"{a}hler $h$ is locally
a product metric where the factors have constant holomorphic sectional curvature -1 (for base curve) or 1 (for the fiber). Following \cite{ADM}, let $(z_{1},z_{2})$ be holomorphic coordinates, and K\"{a}hler form is
\begin{align}
\omega_{h}=2\sqrt{-1}\frac{dz^{1}\wedge d\overline{z}_{1}}{(1-|z_{1}|^{2})^{2}}+2\sqrt{-1}
\frac{dz^{2}\wedge d\overline{z}_{2}}{(1+|z_{2}|^{2})^{2}}=\omega_{1}+\omega_{2}
\end{align}
Let $e$ is unitary frame, such that $e_{i}$ is parallel to $\frac{\partial}{\partial z_{i}}$. Then under frame $e$, we have
\begin{align}
R_{1\overline{1}1\overline{1}}=-1;\,\,\,\,R_{2\overline{2}2\overline{2}}=1;
\end{align}
and others curvature is zero. So the Ricci curvature $R_{11}=-1$, $R_{2\overline{2}}=1$, others is zero. Denote $R_{i\overline{j}}=(-1)^{i}\delta_{ij}$. Since $g$ has  constant 2th-mixed curvature $\mathcal{C}^{(2)}_{\alpha,\beta}=c$. By Lemma \ref{lem3.2}, we can deduce that
\begin{align}\label{3.40}
-(\alpha+\beta)-2\beta F_{1\overline{1}}=2\alpha\Delta F+ce^{2F};\,\,\,\, (\alpha+\beta)-2\beta F_{2\overline{2}}=2\alpha\Delta F+ce^{2F};\,\,\,\,\beta F_{1\overline{2}}=0
\end{align}
Since $\beta\neq0$, so $F_{1\overline{2}}=0$. Setting $F_{1\overline{1}}=x$, $F_{2\overline{2}}=y$. So $y-x=\frac{\alpha+\beta}{\beta}$, $x_{\overline{2}}=F_{1\overline{1}\overline{2}}=0$ and $y_{\overline{1}}=F_{2\overline{2}\overline{1}}=0$. Then we have
\begin{align}
\Delta x=x_{1\overline{1}}=y_{1\overline{1}}=0;\,\,\,\, \Delta y=0
\end{align}
It implies that $\Delta F$ is constant. Hence $F$ is constant and $g$ is K\"{a}hler. Note that at this case we can only have  $c=0$ and $\alpha+\beta=0$.

For case (3), by Lemma \ref{lem3.4},  the Hermitian metric $g$ does not exist on isosceles Hopf surface. This completes the proof of Theorem \ref{th1.2}. \qed

\vspace{0.2cm}

\subsection{High-dimensional Hermitian manifold}

Finally, we present the proof of Theorem \ref{th1.3}.

\noindent{\bf{{\em Proof of Theorem \ref{th1.3}.}}} Since $g$ satisfies $\mathcal{C}^{(k)}_{\alpha,\beta}=0$. It follows from the average trick that
\begin{align*}
&\frac{1}{vol(\mathbb{S}^{2n-1})}\int_{ Z \in T^{1,0}_{p}M, | Z | =1 }\alpha Ric^{(k)}(Z,\overline{Z})+\beta H(Z) d\theta(Z) \nonumber \\
&
\frac{\alpha}{n}u+\frac{\beta}{n(n+1)}(u+v)=\frac{((n+1)\alpha+\beta)u+\beta v}{n(n+1)}=0,\,\,\,\,\,\,\,k=1,2\nonumber
\end{align*}
and
\begin{align*}
&\frac{1}{vol(\mathbb{S}^{2n-1})}\int_{ Z \in T^{1,0}_{p}M, | Z | =1 }\alpha Ric^{(k)}(Z,\overline{Z})+\beta H(Z) d\theta(Z) \nonumber \\
&
\frac{\alpha}{n}v+\frac{\beta}{n(n+1)}(u+v)=\frac{((n+1)\alpha+\beta)v+\beta u}{n(n+1)}=0,\,\,\,\,\,\,\,\,k=3,4\nonumber
\end{align*}
If $\alpha(n+1)+2\beta=0$, then we have $u-v=0$ at every point. And By
\begin{align}
\int_{M}(u-v)d\nu=\int_{M}|\eta|^{2}d\nu \nonumber
\end{align}
we have $\eta=0$, so $g$ is balanced metric. \qed

\vspace{0.5cm}

\noindent\textbf{Acknowledgement.} The authors would like to thank Professors Lei Ni and Fangyang Zheng for many helpful discussions.

\end{document}